# Partial State-Feedback Reduced-Order Switching Predictive Models for Next-Generation Optical Lithography Systems

Raaja Ganapathy Subramanian[a,b,1,*], Barry Moest[b,2], Bart Paarhuis[b,3]

[a] *Eindhoven University of Technology, The Netherlands.*
[b] *ASML NV, The Netherlands.*

---

**Abstract**

This paper presents a partial state-feedback reduced-order switching predictive model designed to support the next-generation lithography roadmap. The proposed approach addresses the trade-off between increasing the number of measurements to improve overlay accuracy and the resulting challenges, including higher measurement noise, reduced throughput and overlay/placement errors under uncertain operating conditions.. By minimizing (die-) placement errors and reducing unnecessary measurements, the method enhances system performance and throughput. This solution employs a streamlined model with adaptive switching logic to manage time-varying uncertainties induced by fluctuating operating conditions. The methodology is implemented on a state-of-the-art lithographic scanner to mitigate the spatial-temporal dynamics of reticle heating, serving as a representative industrial application. Reticle heating, which worsens with increased throughput, introduces spatial-temporal distortions that directly degrade die placement accuracy. Experimental results demonstrate significant improvements: placement errors are reduced by a factor of $2-3$x, and throughput is improved by 0.3seconds per wafer. Importantly, the method accounts for the fact that increased throughput can exacerbate reticle heating, which directly impacts overlay performance. By actively compensating for these thermomechanical effects, the proposed approach ensures that overlay accuracy is maintained or improved – even under increased throughput conditions – highlighting its potential for broader application in advanced lithographic systems, particularly in thermal and vibration control.

---

*Corresponding author
*Email addresses:* raaja.ganapathysubramanian@asml.com, r.g.subramanian@tue.nl (Raaja Ganapathy Subramanian), barry.moest@asml.com (Barry Moest), bart.paarhuis@asml.com (Bart Paarhuis)
[1] Mobile number: +31 682613549

## 1. INTRODUCTION

Current era is dominated by information technology, which is widely considered to have truly begun with the creation of transistors followed by the development of an integrated circuit (IC). Since then, information technology has experienced rapid exponential growth. Similarly, electronics has transitioned from using circuits made of discrete components soldered to a printed circuit board to ICs comprising billions of interconnected transistors on a single chip. These chips provide alternatives to traditionally large and costly discrete optical devices and allow for new applications that traditional techniques could not achieve. This expansion has been driven by exponential advances in computing power, data storage, and communication. This remarkable growth is commonly referred to as Moore's law [1, 2], despite not being a physical law and serving more as an economic hypothesis. For more than 50 years, Moore's prediction has acted as a guiding principle for the semiconductor industry.

Photolithography is the most essential technique in the entire semiconductor production process for determining the scale at which chip elements can be created. It is responsible for transferring the layout of the circuit to the raw material, eventually transforming it into ICs using sophisticated machines [3], as illustrated in Fig. 1. In photolithography, the illuminator generates a (deep/extreme) ultraviolet light (DUV/EUV) beam, which is manipulated through an optical system to interact with a pattern on a quartz/zerodur plate, known as the reticle, thus forming an image on a thin layer of photosensitive material (resist) atop a circular silicon substrate, the wafer. Typically, the minimum feature size in this projection ranges from 3 [nm] to 500 [nm], achieving precision in the order of subnanometers, commonly termed overlay performance or placement errors [4]. Economically, high throughput is essential [5, 6, 7, 8], requiring high operating speeds that introduce disturbances that affect desired performance. To push the limits in the lithography industry, a shift toward model-based prediction and control has been implemented to achieve high-performance goals using raw data [9, 10, 11, 12]. However, these systems can experience various disturbances that cause non-linear behavior. The limitations of linear prediction and control, affected by the "waterbed effect [13, 14, 15, 16]," have prompted a movement towards nonlinear methods, including nonlinear variable gain control [16], filters [15, 11, 10], and switched controllers [13, 14]. While these methods have improved time-domain performance, they often neglect spatial-temporal dynamics that are critical in lithography processes [11]. Despite research in this field [17, 18], practical application of advanced models and control is limited due to challenges such as system complexity and unpredictability as mentioned below.

[2]Mobile number: +31 614611085

[3]Mobile number: +31 623016707

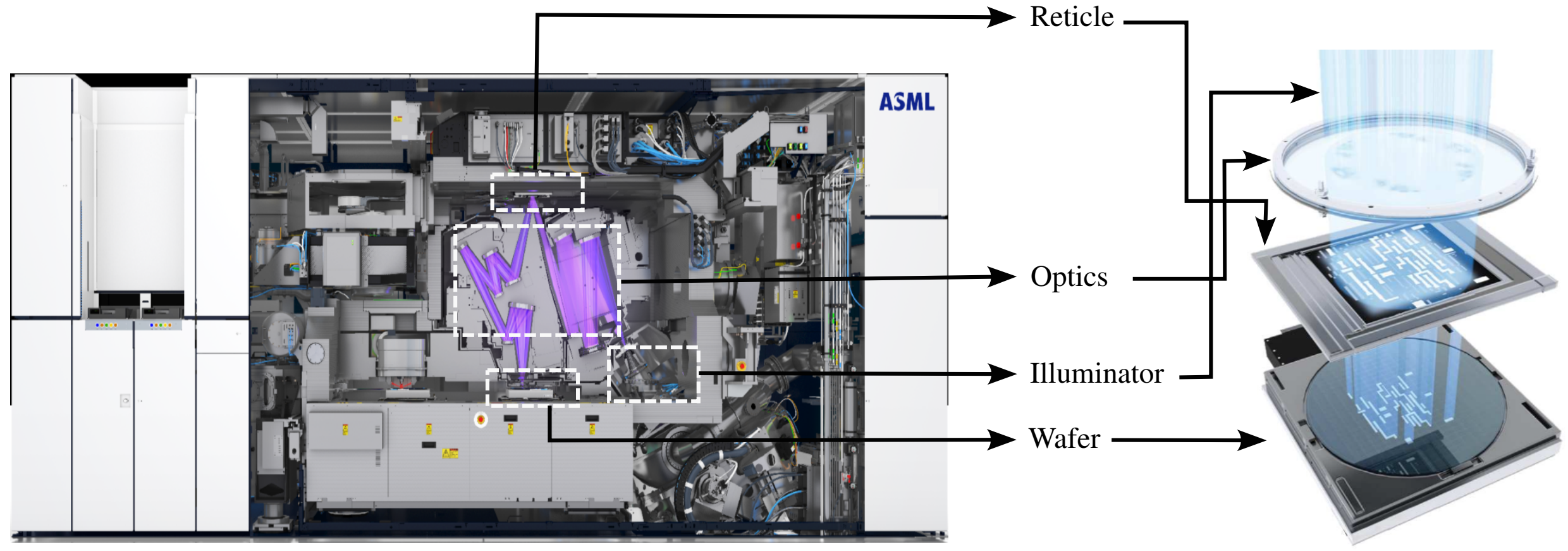

Figure 1: An illustration showing an ASML Twinscan and the exposure process with a lens, reticle, and wafer from top to bottom [3].

1. Mismatch between simulation assumptions and real-time system behavior, such as the use of full versus partial state feedback.

2. Sensitivity to uncertainties and external disturbances, which can degrade control performance.

3. Difficulty in ensuring global closed-loop stability during real-time operation.

4. Trade-offs between accurate physical modeling and model simplicity, which affect usability and interpretability.

5. Concerns related to modularity, scalability, and implementation cost, which impact the feasibility of deployment.

The first three factors directly influence system performance, specifically concerning overlay and placement errors, while the fourth and fifth factors serve as practical limitations that often restrict the feasibility of suggested methods. Simply put, the objective of this work is to attenuate placement errors arising from these spatial-temporal behaviors under a typical production situation (i.e., under uncertainties), while simultaneously enhancing system productivity.

Further to bridge the gap between theoretical control concepts and practical implementation in lithography systems – while maintaining both throughput and overlay performance – we propose a partial state-feedback, reduced-order switching predictive control approach. This method is designed to effectively manage system uncertainties, enhance performance, and remain globally stable, computationally efficient, and intuitively implementable. The main contributions of this work are as follows:

- To address challenges 1, 4, and 5, the use of reduced-order linear time-invariant (LTI) models derived from nominal large-scale system dynamics (in this case a finite element model) is proposed. By selectively decoupling only the essential spatial-temporal dynamics, computational efficiency without compromising model fidelity is ensured.
- To address challenge 2, a switching logic mechanism is introduced that adapts to uncertainties arising from physical effects, thereby improving robustness.
- To address challenge 3, the proposed control strategy is proven to ensure global uniform ultimate bounded asymptotic stability (GUAS), supported by the small-gain theorem for nonlinear systems.
- To validate the proposed approach, it is implemented on a state-of-the-art lithography system, with a specific focus on mitigating spatial-temporal distortions induced by reticle heating.

The structure of this paper is organized as follows. Section 2 presents the problem formulation and connects it to the key challenges outlined in the Introduction regarding spatial-temporal dynamics in lithography systems. Section 3 introduces the proposed control methodology and analyzes its stability properties. Section 4 and 5 demonstrates the effectiveness of the proposed approach through experimental validation on a state-of-the-art lithography platform. Section 6 concludes the paper with a summary of key findings, and Section 7 outlines potential directions for future research.

## 2. Problem Formulation

As discussed in the Introduction, spatial-temporal variations in optical lithography systems can lead to sub-nanometer misalignments, which directly degrade patterning precision. While measuring these deformations during wafer exposure can help mitigate their impact, doing so often results in reduced productivity. A critical example of such deformation is reticle heating. Since the reticle contains the mask that defines the chip layout, any thermal-induced deformation of the reticle directly affects the accuracy of the pattern transferred onto the wafer, thereby impacting the overall integrity of the integrated circuit (IC) production process.

During standard production, light is absorbed by the reticle as it progresses toward the lens, as shown in Fig. 2. The reticle is usually positioned on a clamp. The chrome absorber layer, which delineates the pattern on the reticle, absorbs this light during exposure, causing it to heat up; this phenomenon is termed reticle heating (RH). Moreover, depending on the specific layer that is exposed, the reticle may feature a pellicle, serving as a protective barrier against contamination. In addition, heat dissipation occurs from both the upper and lower sides of the reticle to the surroundings, and there is a cooling flow designed to mitigate the effects of reticle heating. In typical manufacturing cycles, a wafer batch is exposed using a reticle, commonly referred to as a lot. Consequently, the reticle undergoes frequent thermal cycling, as depicted in Fig. 3. During this process, the temperature of the reticle rises upon exposure of each wafer (approximately tens of seconds), and it does not sufficiently cool during the wafer exchange period (less than a few seconds). This results in a significant temperature gradient between the first and last wafers in a lot. Overall, the global heating trend follows an exponential progression over time, leading to distortions in the image projected onto the wafer. Additionally, once a lot is fully processed, not only do the wafers need replacement with a new batch, but the reticles do as well. This replacement involves unclamping and reclamping the reticle from the reticle stage, introducing variability that alters the thermomechanical boundary conditions of reticle heating dynamics (RH). As shown in Fig. 3, RH physics also demonstrates spatial-temporal characteristics [19, 20, 21, 22]. The extent of these dominant spatial-temporal distortions presents a significant challenge, particularly for high-throughput and high-performance systems. The primary focus will be on placement/overlay errors in horizontal direction and the vertical direction will be ignored for simplicity.

### *2.1. Control Objective*

Having discussed the background of the physical system, in this subsection the control objective will be discussed. Let $x \in \mathbb{R}^n$ represent the system state (e.g., the predicted full reticle deformation, $u_e \in \mathbb{R}^m$ denote the exogenous input (e.g., heat source), and $z \in \mathbb{R}^p$ be the output associated with placement error. The control objective is to minimize $z$, such that

1. the predicted misalignment remains bounded under uncertainties:

$$||z|| \leq \epsilon, \quad \text{(e.g., sub-nanometer).} \tag{1}$$

   The measurements available, $y \in \mathbb{R}^q$, represent sparse sensor data collected (as shown in Fig. 4) during alignment and the associated feedback to the model as $u_f$, typically available only at start of each wafer. This partial feedback is used to estimate and update $x$.

2. measuring each alignment marks can take up to 0.3 seconds. Therefore, an additional objective is to reduce measurement time by selectively skipping measurements on selected marker sets out of the full marker set in as shown in Fig. 4 — this would include edge markers but also includes skipping top and/or bottom row. Further, using the models to predict the deformation based on partial feedback, while simultaneously improving placement accuracy.

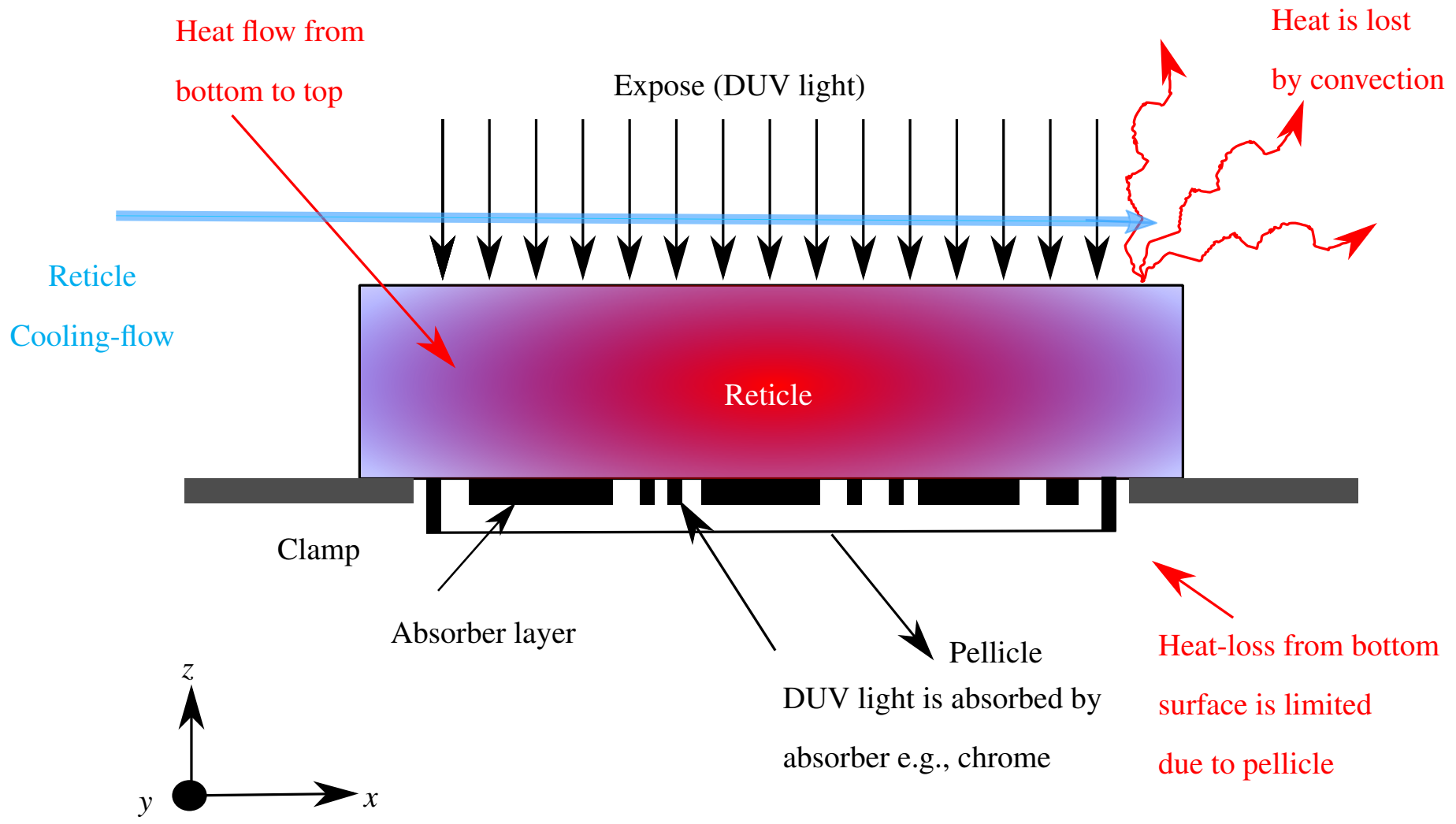

Figure 2: Artistic impression of reticle heating physics.

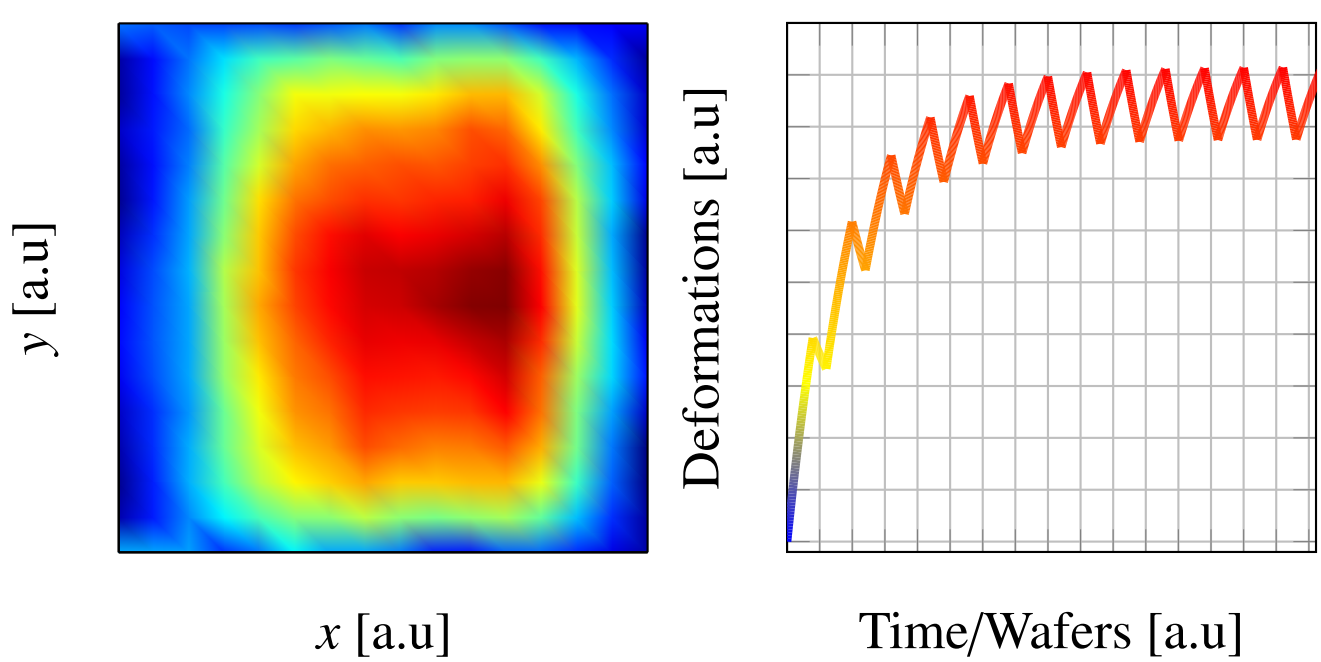

Figure 3: Spatial-temporal dynamics of the reticle heating.

## 2.2. *State-space formulation*

Now that the sparse measurements from alignment sensors are taken in feedback and that the model has to deal with uncertainties, the state-space is formulated as follows:

$$\hat{\mathcal{P}} := \begin{cases} \dot{x}_i = A_i x_i + B_i u_e + B_f u_f, \\ u_f = \Gamma(.) y, \\ u_i = C_i x_i. \end{cases} \tag{2}$$

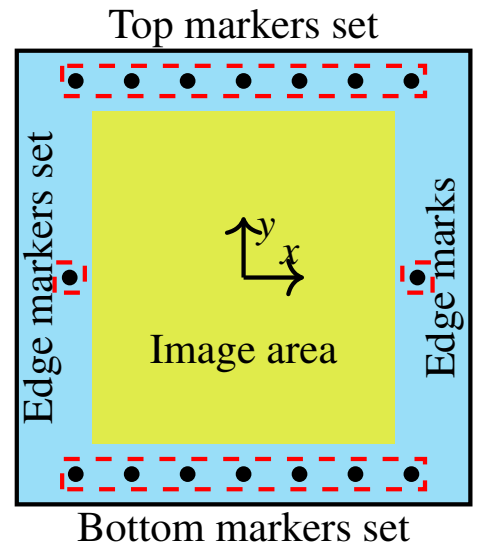

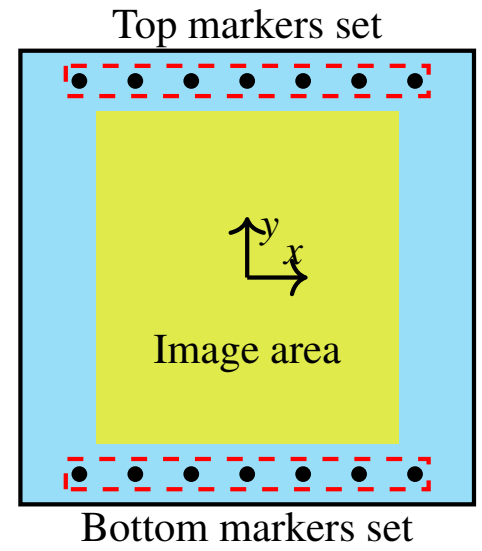

Figure 4: Exposure layout with varying alignment marks; reticle is exposed in the image area.

where the index $i$ depends on the selected model per uncertainty classification. Associated variables are defined as:

- $x_i$: internal deformation state
- $u_e$: external heating profile
- $u_f$: feedback signal from alignment sensor
- $z$: actual placement error
- $\Gamma(.)$: converts sparse feedback measurements into an appropriate dense layout in the same density as $z$, using application-specific inputs (like image-area etc.)
- $y$: sensor measurements available at runtime (subset of $z$)
- $u_i$: predicted placement error that acts as a correction signal mimicking the physics (in this case the reticle heating)
- $\Delta$: is the uncertainty that occurs in a real-time environment (due to changes in boundary /operating conditions etc.)
- $C_i$: represents the output mapping from $x_i$ to $u_i$. In the context of reticle heating application, this mapping serves as a transformation from temperature to deformation, thereby, comprising the mechanical boundary conditions.

and a description of such a system as shown in Fig. 5 falls under the generalised plant framework [23].

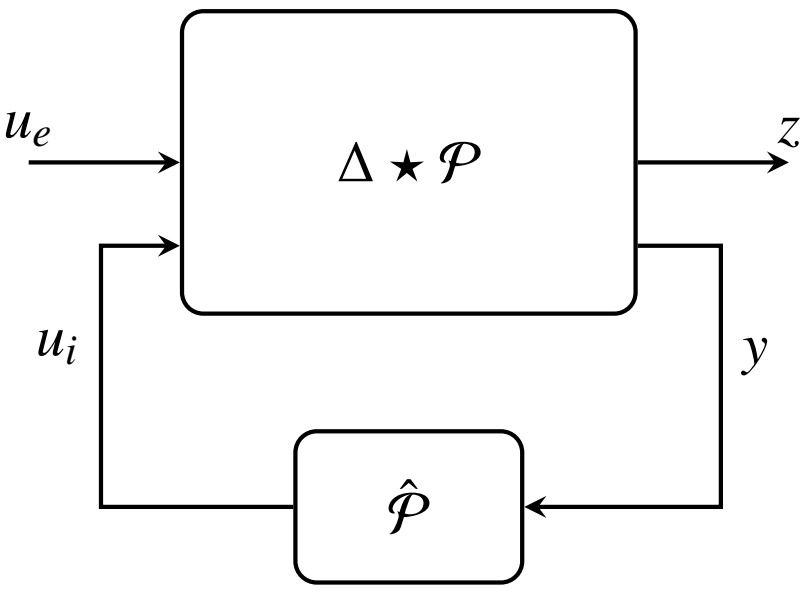

Figure 5: Simplified depiction of a partial state-feedback model in a generalized plant framework. Note that the $\star$ denotes the LFT operator.

## 3. DESIGN OF PARTIAL STATE-FEEDBACK REDUCED-ORDER SWITCHING MODELS

To address the problem of spatial-temporal behavior introduced in Section 2, this section begins with a step-by-step design philosophy in Section 3.1 and covers the closed-loop dynamics stability analysis in Section 3.2.

### *3.1. STRUCTURE OF THE PROPOSED METHODOLOGY*

To implement the proposed methodology in real time, four main steps are required:

*Step 1 – Integrating the Scheduler* $\Phi$ *and Closing the Feedback Loop.* Examine the closed-loop control configuration depicted in Fig. 6, where the plant $\mathcal{P}$ functions in response to an external input $u_e$. This input is affected by a time-dependent uncertainty $\Delta$ through an uncertain input channel $u_\Delta$. The setup includes a reduced-order predictive switching model $\hat{\mathcal{P}}$ that uses partial state-feedback based on measured events $y$ and aims at the desired result $z$. The predictive model also incorporates the scheduler $\Phi$, which employs the history describing current and previous events, represented as $\mathscr{I} \in \{y, u_e, u_\Delta\}\ \forall T \in \{0, \ldots, t_n\}$, to alternate between the applicable models. The primary focus is on history-based switching, which selects the appropriate models using available information from $y, u_e, u_\Delta$, as opposed to variable-gain control based that only changes the gain solely based on the magnitude of the feedback signal [24, 25, 26, 16, 27] or its product with the signal and its derivative [28, 29, 30]. The information-based switching strategy is used because of its event driven nature [31, 32]. Adjusting the control based on other means would complicate the design. Instead, historical data $\mathscr{I}$ is used which is always available and its often inherent to physics of spatial-temporal system, to guide model selection efficiently. In order to incorporate this in the proposed scheme, the scheduler $\Phi$ is chosen as follows:

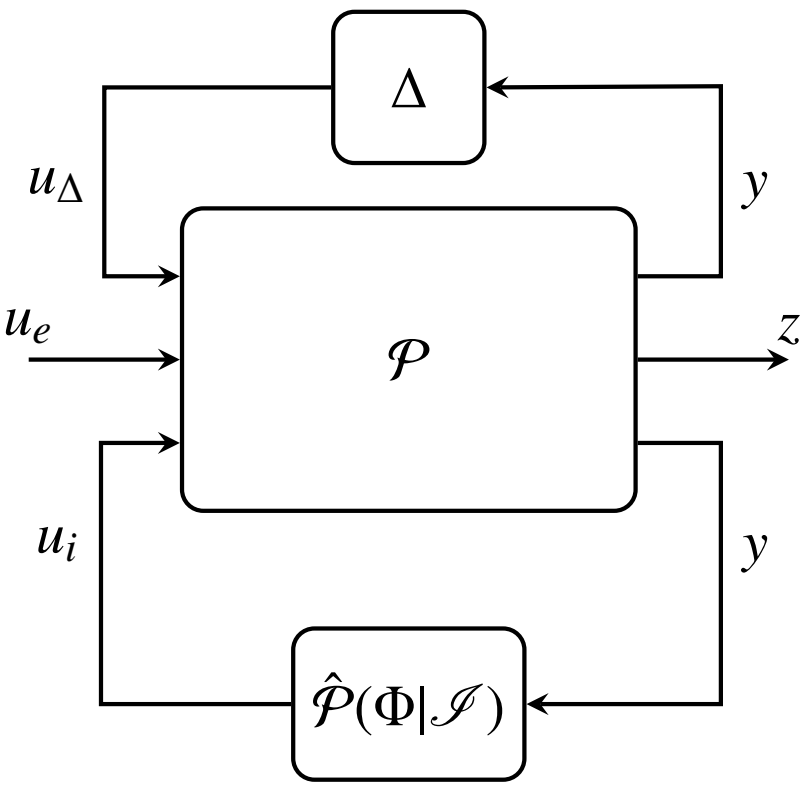

Figure 6: Simplified depiction of a partial state-feedback model with order reduction and switching in a closed-loop system.

$$\Phi = \{\mathcal{M}_i|\mathscr{I} \in \mathscr{R}_i \forall i \in \mathbb{Z}\} \tag{3}$$

*Definition* 3.1. An operating condition, referred to as a regime $\mathscr{R}$, is defined by the integration of a feedback signal $y$, an external signal $u_e$, and uncertainty $u_\Delta$. The regime $\mathscr{R}_i$ for $i = 0$ is identified as the nominal regime and is associated with a nominal model $\mathcal{M}_0$. Conversely, for $i \neq 0$, it is described as the uncertain regime, which involves a set of uncertain models $\mathcal{M}_{i\in\mathbb{Z}\setminus\{0\}}$.

Consider an example using (3) with two distinct regimes: the nominal regime $\mathscr{R}_0$ and its associated model $\mathcal{M}_0$, and another regime $\mathscr{R}_1$, identified through $\mathscr{I}$. In this scenario, within regime $\mathscr{R}_1$, the set $\{u_e, u_\Delta\} \in \emptyset$. Initially, with model $\mathcal{M}_0$, the target performance $z$ is achieved optimally when the system operates within regime $\mathscr{R}_0$. However, if the system transitions to regime $\mathscr{R}_1$ due to the impact of $y$, the optimal performance of $z$ is not maintained unless the scheduler $\Phi$ engages model $\mathcal{M}_1$, appropriate for regime $\mathscr{R}_1$. It is important to note that the choice of the nominal model greatly affects the performance of $z$. Thus, by effectively configuring the scheduler $\Phi$ and selecting the nominal model $\mathcal{M}_0$, the advantages are combined to ensure optimal nominal and robust performance.

*Step 2 – Formulation of Nominal Reduced-Order Model.* Partial differential equations govern spatial-temporal systems, as noted in [33, 18, 34], and are typically derived using finite-element models (FEM) or computational fluid dynamics (CFD) models. These large-scale models are often too complex for real-time application due to their high dimensionality (often exceeding $1e^6$) necessary for precise accuracy [11]. As an illustration, one method detailed in [35] runs a dynamic FEM alongside the physical system, providing real-time forecasts of the nominal spatial-temporal dynamics. This approach presents numerous modeling challenges, not only from theoretical and experimental standpoints, but also in ensuring robust performance.

To address this issue, model-order reduction techniques have been introduced, such as Balanced Truncation [36], Krylov's Subspace Methods [37], and the Proper Orthogonal Decomposition (POD) method [38]. These techniques are prevalent in various engineering domains. However, even though they decrease the dimensions of the FEM / CFD models, they often lose physical interpretability and continue to lack computational efficiency [18].

To achieve this, a parametric reduction is applied leveraging a Krylov's moment-matching method [39] that maintains the original physical interpretation of the reduced-order model. This approach can be expressed as ,

$$\mathcal{M}_i := \begin{cases} \dot{x}_i = A_i x_i + B_i u_e, \\ u_i = C_i x_i. \end{cases} \tag{4}$$

Here, $u_i$ in (4) is the estimated behavior of the underlying process, as shown in Fig. 6. The aim is to ensure that $||z|| \to 0$ and its also important to observe that in equation (3), $\mathcal{M}$ denotes the reduced-order model corresponding to different regimes $\mathscr{R}$ of operating conditions. Upon the scheduler's detection of a regime shift, the model is updated, and the corresponding internal states are transferred from $\mathcal{M}_k$ to $\mathcal{M}_{k+1}$.

*Step 3 – Using the partial feedback signal y to ensure that the model $\mathcal{M}_i$ is up-to-date.* In industrial settings, complete state feedback often demands considerable time. Consequently, what is typically available is partial feedback, characterized by being limited and generally representing only a subset of the overall layout of physical effects. As a result, models deployed in such scenarios must accommodate these constraints. To integrate partial state-feedback into the reduced-order models discussed in the previous section, (4) is enhanced as in equation (2).

This is achieved by explicitly incorporating an input channel $B_f$ for feedback, with $y$ representing the actual feedback measurements. In addition, in (2), a mapping function $\Gamma(.)$ is introduced. This function is designed to convert sparse feedback measurements into a suitable dense layout, using application-specific inputs. It is important to note that the variable $u_i$ in (4) and (2) is intentionally defined as such. The key lies in the definition of the information set $\mathscr{I}$, which is passed into the functional $\hat{\mathcal{P}}(\Phi|\mathscr{I})$. This set comprises the exogenous input $u_e$, the measured output $y$, and the uncertainty term $u_\Delta$. To reformulate the problem in the 9-gang/generalized plant framework [23], the partial output measurement $y$ is incorporated into the feedback loop, and the information set $\mathscr{I}$ is treated as a variable within the functional representation of $\hat{\mathcal{P}}$.

*Step 4 – Transforming reduced-order models $\mathcal{M}_i$ to incorporate uncertainty.* To efficiently manage uncertainties as highlighted in Section 3.1, the framework is divided into a nominal model, denoted as $\mathcal{M}_0$, alongside a collection of uncertain models, represented by $\mathcal{M}_\Delta$. By employing a centering method [23], these models are expressed in the form:

$$\hat{\mathcal{P}}(\Phi|\mathscr{I}) = \mathcal{M}_n + \mathcal{M}_\Delta = \mathcal{M}_n + \Delta_m(\Phi|\mathscr{I}), \quad \text{where } m = \arg\max_m \Phi. \tag{5}$$

Here, $\mathcal{M}_n$ serves as the nominal model, while $\Delta_m(.)$ represents the set of uncertain dynamical models. It is apparent from the equation (5) that the collection of models belongs to the category of Lur'e-type systems [40], which integrate both a linear dynamical component and a nonlinear/uncertain one (a Lur'e-type system description).

Furthermore, this approach ensures that the uncertain dynamics is bounded, specifically $||\Delta_m(.)||_\infty \in \mathbb{R}$, and modulated by the scheduler $\Phi$ according to (3) [31, 32]. Finally, the scheduler mechanism can explicitly defined as follows:

$$\Phi := \begin{cases} i \in \mathbb{Z} \\ \text{if } \mathscr{R}_i(\mathscr{I}_t) \geq \mathscr{R}_j(\mathscr{I}_{t^-}),\ \forall j \in \mathbb{Z}_{\setminus\{0\}}, t \in \mathbb{R}_{\geq 0} \\ \mathscr{I}(.) \in \{y, u_e, u_\Delta\}. \end{cases} \tag{6}$$

Considering the reticle heating application in lithography system, the reticle transitions between various subsystems — such as the exposure stage, metrology station, and reticle handler. Each of these locations introduces different thermal and mechanical conditions due to variations in airflow, illumination, and mechanical contact. These transitions affect the system's behavior and, conse-quently, the control strategy. In this context, the regime $\mathscr{R}$ can be defined based on the historical and real-time information $\mathscr{I} \in \{y, u_e, u_\Delta\}$. Here, $y$ is a measured artefacts and what is known from sparse measurements, $u_e$ could be environmental or equipment-related inputs (e.g., stage velocity, chuck temperature), and finally $u_\Delta$ can be detected changes in system configuration or operating mode. For example, when the reticle is moved from the exposure stage to a non-exposure zone (e.g., for inspection or handling), the scheduler $\Phi$ may detect a shift in $\mathscr{I}(.)$ and classify the system as entering a new regime $\mathscr{R}_n$, prompting a switch to a different model or control strategy. This allows the system to adapt to location-specific boundary conditions without requiring direct measurement at every step [31, 32]. Therefore, The specific definition of the regime $\mathscr{R}$ is application-dependent and will vary based on the context, objectives, and characteristics of the input signals $\mathscr{I}(.) \in \{y, u_e, u_\Delta\}$. As such, $\mathscr{R}$ should be designed to reflect the performance criteria or priorities relevant to the particular system under consideration.

### 3.2. STABILITY ANALYSIS

To conduct a stability analysis of the dynamic behavior influenced by the proposed approach, Fig. 6, (2) and (5) is utilized, along with the generalized-plant framework as defined in [23]. The the closed-loop system is transformed in to a Lur'e-type system, illustrated in Fig. 7. The closed-loop dynamics can be minimally realized by integrating the uncertainty with the nominal predictive model, employing upper and lower linear-fractional-transformations (LFT) as outlined respectively in [23]:

$$\begin{bmatrix} \dot{x} \\ y \\ z \end{bmatrix} = \begin{bmatrix} \mathcal{A} & \mathcal{B}_e & \mathcal{B}_f & \mathcal{B}_i \\ \mathcal{C}^y & \mathbb{O} & \mathcal{D}_f^y & \mathcal{D}_i^y \\ \mathcal{C}^z & \mathbb{O} & \mathcal{D}_f^z & \mathcal{D}_i^z \end{bmatrix} \begin{bmatrix} x \\ u_e \\ u_f \\ u_i \end{bmatrix}, \tag{7a}$$

$$u_i = \Delta_m(\Phi|\mathscr{I}), \tag{7b}$$

with the uncertainty defined as:

$$\Delta := \begin{cases} \dot{x}_\Delta = A_\Delta x + B_\Delta u_f, \\ u_f = \Gamma(.)y, \\ u_\Delta = C_\Delta x_\Delta + D_\Delta u_f. \end{cases} \tag{8}$$

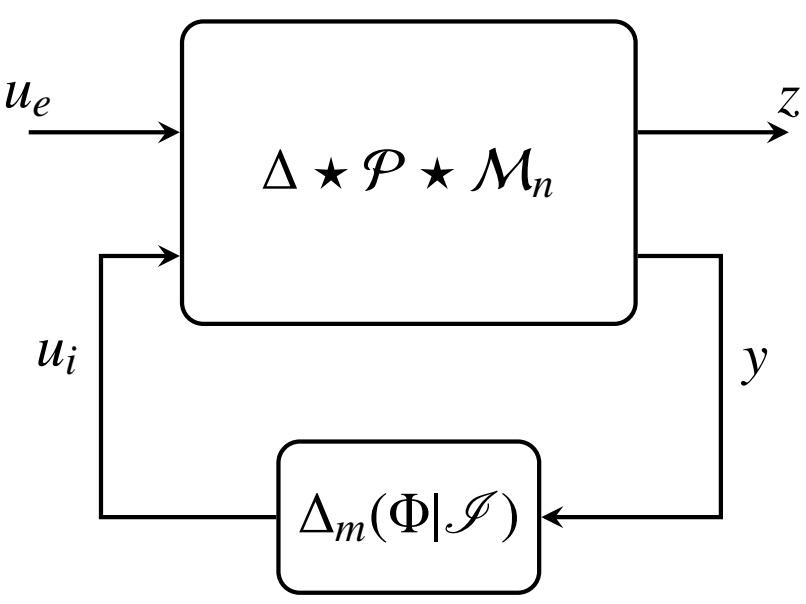

Figure 7: Schematics of the proposed methodology in a Lur'e-type system.

*Assumption* 3.1. The closed-loop dynamics of the actual plant, when compensated using the nominal reduced-order partial state-feedback denoted as $\Delta \star \mathcal{P} \star \mathcal{M}_n$, fits within the category of globally uniformly ultimately bounded asymptotically stable systems.

Let us consider $x^\star$ to be the equilibrium point of (7) such that $||z|| \rightarrow 0$. It is important to note that $x^\star$ is the sole equilibrium point with the property that $||z||= 0$, because the minimum state-space representation in (7) ensures observability. This means that the observability matrix is of full rank, and hence, the system of equations $z = 0, \partial z/\partial t = 0, \ldots, (\partial^{n-1}z)/(\partial t^{n-1}) = 0$ yields a unique solution $x^\star$, when $||z|| \rightarrow 0$. Using the small-gain theorem as specified in [13, 24, 41] together with the stability theorem presented in [16], the global uniform ultimate bounded asymptotic stability can be guaranteed. As a result, under such conditions, it is possible to precisely forecast the dynamics of the plant, which implies that $||z||$ will be reduced to zero for all $t \in \mathbb{R}$.

## 4. LITHOGRAPHY SYSTEM APPLICATION: RETICLE HEATING

In this section, the proposed methodology will be implemented on a high-precision industrial optical system, specifically an ASML Twinscan, as shown in Fig. 1 for the reticle heating application as shown in Fig. 2 in the problem formulation section (see, Section 2). Before proceeding to discuss the experimental results in Sections 4.2, first the stability conditions are established in Section 3.2, as detailed in Section 4.1.

### *4.1. STABILITY CONDITIONS*

The reticle is constructed from physical materials and during its heating process, energy absorption results in deformations. Thus, reticle heating displaying unbounded dynamics is practically unlikely. Applying similar reasoning, analogous conclusions for the nominal model $\mathcal{M}_n$ is drawn and the specific models that represent the uncertain dynamics characterized by $\Delta_m(\Phi|\mathscr{I})$. Furthermore, utilizing the design methodology described in Section 3.1, 3 moments are identified that describe both the nominal and uncertain reticle heating models. Consequently, it is reliable to say that the dynamics of the actual system and the nominal model in a closed loop, namely, $\Delta \star \mathcal{P} \star \mathcal{M}_n$, comply with the assumption 3.1. In addition, the uncertain model definition given by (7) is also adopted. Therefore, as demonstrated in Section 3.2, all necessary conditions are met to ensure that the equilibrium point $x^\star$ (where $||z|| = 0$), reached by applying the proposed method (7), is GUAS. ■

### *4.2. EXPERIMENTAL RESULTS*

To evaluate and measure how the proposed method performs with respect to throughput and overlay placement errors, a series of real-time measurements is obtained, as depicted in Fig. 8. In this experiment, initially, during the regime labeled as $\mathcal{R}_a$, the heating of the reticle shows standard behavior (without uncertainty) and gradually heats up to its saturation level. Subsequently, the reticle used in $\mathcal{R}_a$ is released, transitioning to a new regime $\mathcal{R}_b$, where different dynamics and boundary conditions apply. Consequently, when the reticle returns for a second series of exposures and is re-clamped, it tends to display the dynamics characteristic of $\mathcal{R}_c$ due to the altered boundary conditions resulting from the re-clamping process. Initially, throughput performance is evaluated relative to the current approach illustrated in Fig. 4, which utilizes top-bottom-edge alignment marks as detailed in [17], with measurements conducted for each wafer. Typically, additional alignments incur time costs as a result of the requirement for physical measurements. In this scenario, each measurement takes approximately 0.3s per wafer. By adopting the proposed methodology, a spatial-temporal prediction model to characterize the RH across the entire reticle is implemented. Leveraging this model allows us to bypass the need for edge-alignment marks, substituting measurements with predictive estimates. Specifically, the measurements from top-bottom alignment marks (see Fig. 4) is used and combined with the models discussed in Section 3 to predict $z$ and delineate the RH physics throughout the reticle. Consequently, this approach saves 0.3s per wafer, which can lead to an increase of up to 7 wafers per hour, a significant improvement for throughput.

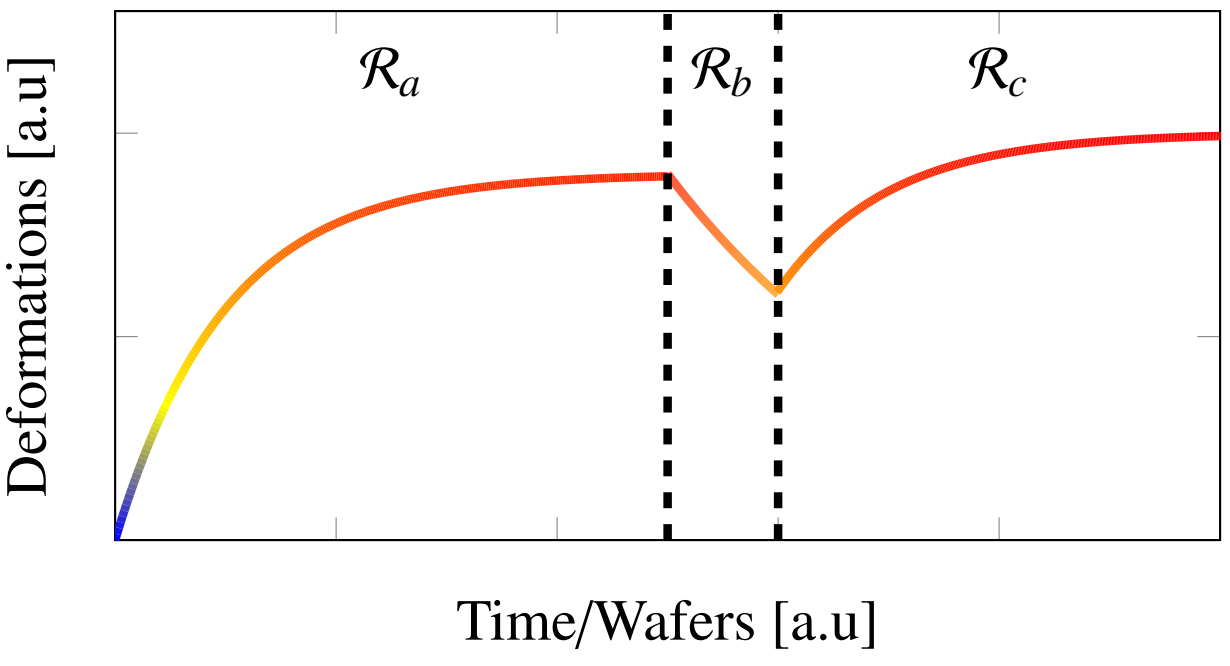

Figure 8: Experiment triggering nominal and uncertain conditions.

In Fig. 9, a comparative analysis of the overlay placement error performance is presented. This analysis is not only conducted against the current standard [17], but also against a linear version of the proposed approach that solely relies on the nominal model $\mathcal{M}_n$ under all conditions. This latter comparison highlights the critical need for careful application of model-based predictions in industrial scenarios, as misuse can lead to significant performance drawbacks. Additionally, an experiment detailed in Fig. 8 was conducted over 2 lots, with each lot comprising of 16 wafers. It is observable that the placement error associated with the proposed methodology exhibits greater stability and improved accuracy in both the $x-y$ direction compared to the current method, under both nominal and uncertain conditions. Furthermore, the linear counterpart of the proposed technique experiences up to 3x reduction in performance due to unknown uncertainties. Notably, with the adoption of the proposed approach, the performance of the first 2 wafers in the $x-y$ direction shows a 2x improve-ment over the conventional method. Additionally, the proposed methodology achieves performance that is either comparable or superior even without the need for measuring additional alignment marks. Consequently, while increasing throughput by up to 7 wafers per hour, the overlay placement errors remain stable and enhanced.

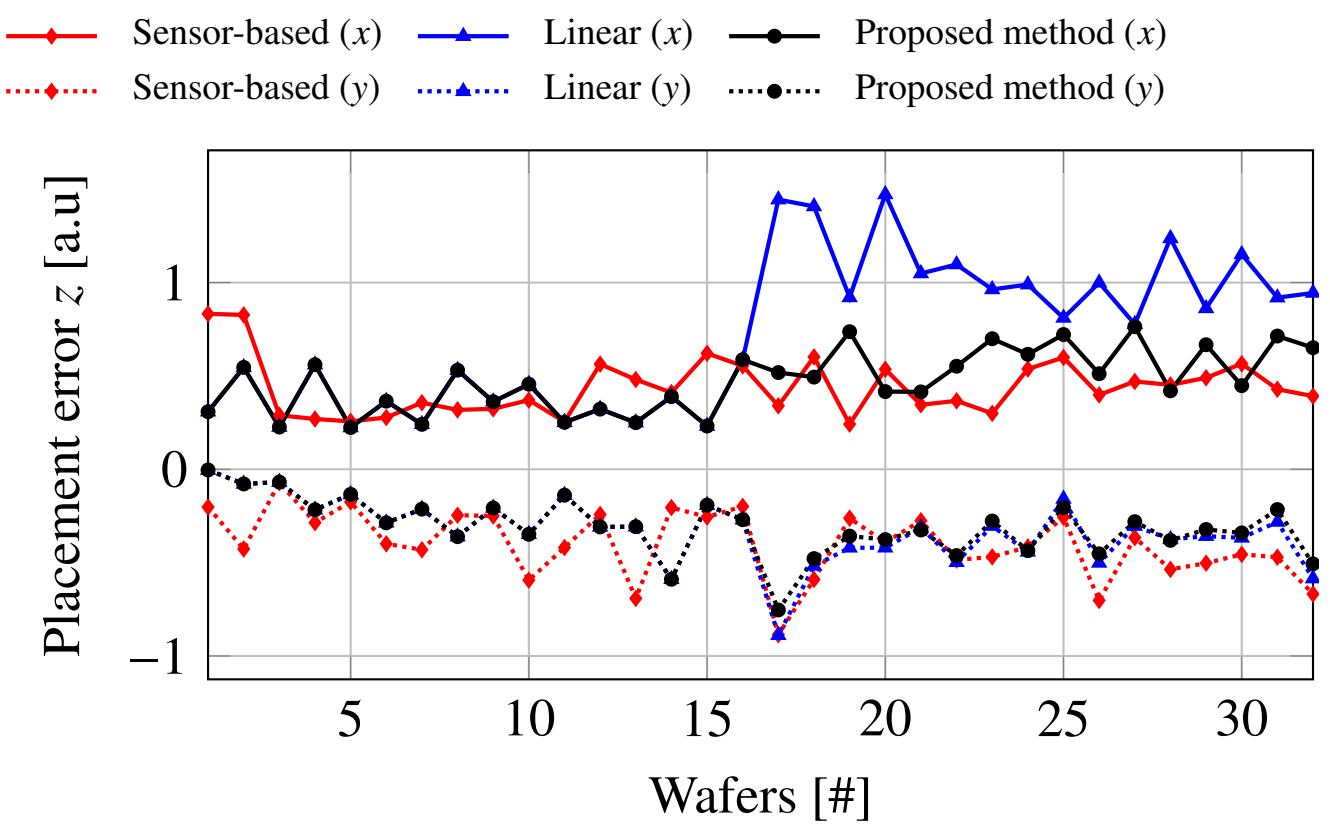

Figure 9: Placement error measurement in $x - y$ to compare the status-quo (red) against the proposed methodology (black) and its linear variant (blue).

*Remark* 4.1. In the context of reticle heating, $\Gamma$ encodes the spatial relationship between alignment mark coordinates and the predicted heating locations, effectively bridging the gap between partial observations and the full-state input required for control. This design enables the model to leverage available feedback while accounting for the spatial sparsity of the measurements. Given this choice, it is good to note that the input channel $B_f$ with respect to $u_f$ reduces essentially to an identity matrix. This reflects a direct mapping between the feedback input $u_f$ and the corresponding state components, which is appropriate given the structure of the reduced-order model and the nature of the partial feedback available in this application (with respect to reticle heating)

*Remark* 4.2. In the proposed switching predictive control framework, $u_\Delta$ is included as part of the switching input set. However, in the specific case of reticle heating, $u_\Delta$ represents a system-generated signal associated with events such as reticle reclamping. This signal is not directly measured; rather, its occurrence is inferred approximately based on system behavior. The timing of such events is only known retrospectively and with limited precision. Despite this, the framework is designed to handle such uncertainty by leveraging historical patterns using other signals from $\mathscr{I}$ and partial observations to infer the active regime, ensuring robustness even when certain inputs like $u_\Delta$ are not explicitly available in real time.

## 5. DISCUSSIONS

In the context of customer interactions with these system during IC production, users may not exclusively utilize the image-area depicted in Fig. 4; they might opt for alternative image-areas as illustrated in Fig. 10. The proposed approach is designed to be robust and remains applicable without the necessity for redesign or retuning despite these on-the-fly layer changes. To demonstrate its efficacy, when conducting the same experiment with a reduced image-area, as illustrated in Fig. 10, Fig. 11 shows that the proposed methodology has improved the performance of both the nominal and uncertain regimes by a factor of 3x, while maintaining consistency. In contrast, similar to the previous findings, the linear version of the proposed method encounters significant challenges, leading to a performance decline of $2 - 3$x compared to both the status quo and the proposed methodology. This demonstrates that the proposed approach allows overlay placement errors to be comparable to or better than existing methods, while simultaneously enhancing system throughput. Additionally, it illustrates how state-of-the-art technologies can be integrated into industrial applications without sacrificing resilience, simplicity, and effectiveness. In conclusion, it is important to emphasize the intuitive nature of the proposed methodology. Linear and uncertain models can be employed via first-principle-based FEM with model reduction and data-driven system identification [42], respectively. Importantly, there are no limitations on the order of the prediction models considered, and only output measurements are utilized. Furthermore, the GUAS of the equilibrium point is ensured and is naturally integrated into the design of the proposed methodology.

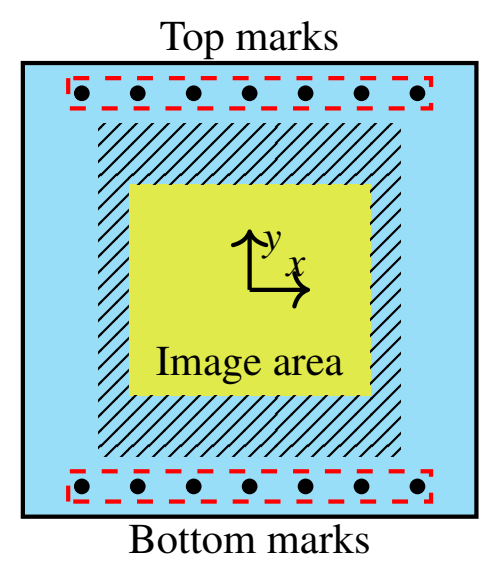

Figure 10: Schematic representation of small-exposure layout.

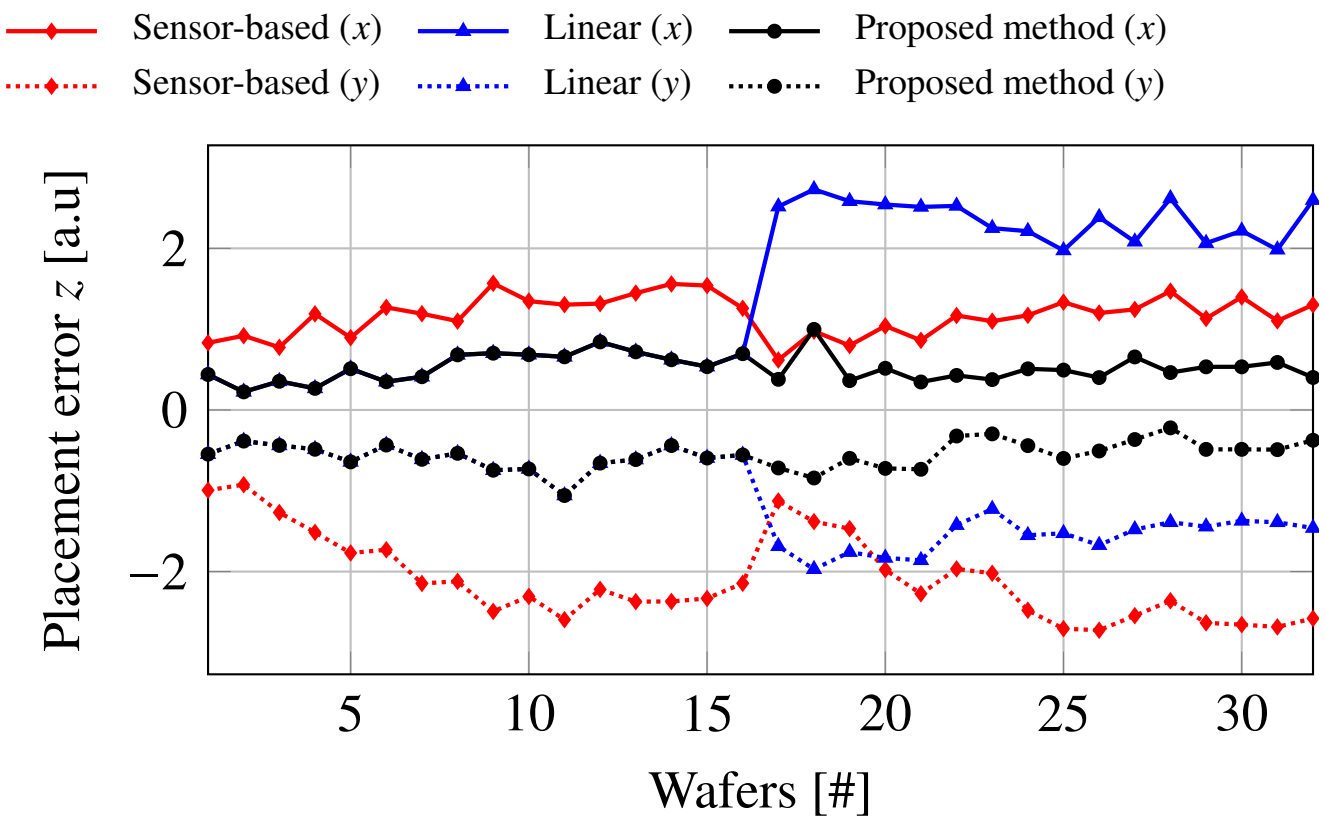

Figure 11: Placement error measurement in $x - y$ for a smaller image area (see, Fig. 10).

## 6. CONCLUSIONS

In this study, a partial state-feedback reduced-order switching predictive model is designed to support the future lithography roadmap. This work addresses the balance between the demand for increased measurements, noise, and overlay errors (both within and across multiple wafers) under uncertain operating conditions. This approach aims to navigate these trade-offs more effectively. This is accomplished by employing reduced-order linear models that can transition among various models based on a scheduling logic derived from historical data, thus handling time-varying uncertainties triggered by operating conditions. To address measurement configurations and maintain practical relevance, a partial state-feedback framework updates the model's internal states using these measurements. Sufficient criteria for global uniform ultimate bounded asymptotic stability through time-regularization techniques is also introduced. A generic execution method automates machine-in-the-loop initialization and execution of the methodology. The strategy is applied to a cutting-edge lithography scanner to mitigate the spatial-temporal dynamics of reticle heating. Experimentally significant improvements are achieved, with placement errors reduced by up to $2 - 3$x and throughput enhanced by 0.3 [s] per wafer in all nominal and uncertain operating conditions. These results are compared to the current standard and the use of only the linear component of the proposed method.

## 7. FUTURE OUTLOOK

As previously mentioned in this manuscript, the relentless demand for high system productivity, coupled with the requirement to maintain equivalent or superior overlay accuracy and placement, influences and imposes demands on the enabling methodologies. These models must, in essence, be resilient and capable of functioning in the midst of uncertainties. There are several key components to consider: 1) A unified predictive framework: although the main focus has been the reticle-heating subsystem, the entire assembly also experiences spatial and temporal variations across both lenses and wafers; 2) achieving enhanced throughput involves identifying and removing unnecessary measurements. Hence, a seamless integration of prediction, learning, and intelligence can allow measurement reductions without compromising overall system efficiency [43, 44, 45]; and lastly, 3) making plausible assumptions is crucial for practical deployment in industrial settings. The greater the disparity, the more challenging it becomes to reconcile theoretical concepts with practical applications. Consequently, the need is effectively summarized as the "development of approaches and techniques for designing adaptable and robust models for high-performance lithographic systems functioning in uncertain conditions while maintaining simplicity." In conclusion, the authors anticipate that the findings reported herein will motivate further exploration and application of nonlinear and/or learning models in real-time systems to enhance performance and bridge the gap between theoretical insights and practical implementation.